\input amstex\documentstyle{amsppt}  
\pagewidth{12.5cm}\pageheight{19cm}\magnification\magstep1
\topmatter
\title {The canonical basis of the quantum adjoint representation}\endtitle
\author G. Lusztig\endauthor
\address{Department of Mathematics, M.I.T., Cambridge, MA 02139}\endaddress
\thanks{Supported by NSF grant DMS-1303060.}\endthanks
\endtopmatter
\document

\define\mpb{\medpagebreak}

\define\frl{\forall}

\define\sqc{\sqcup}

\define\qua{\quad}

\define\op{\oplus}
   
\define\part{\partial}

\define\ra{\rangle}
\define\n{\notin}

\define\m{\mapsto}
\define\do{\dots}
\define\la{\langle}

\define\sub{\subset}    

\define\T{\times}
\define\ti{\tilde}
\define\nl{\newline}
\redefine\i{^{-1}}
\define\fra{\frac}
\define\un{\underline}
\define\ov{\overline}

\define\Hom{\text{\rm Hom}}

\define\a{\alpha}

\redefine\d{\delta}
\define\e{\epsilon}
\define\et{\eta}

\redefine\l{\lambda}

\redefine\L{\Lambda}

\redefine\AA{\bold A}
\define\BB{\bold B}

\define\NN{\bold N}

\define\QQ{\bold Q}

\define\UU{\bold U}

\define\ZZ{\bold Z}

\define\ca{\Cal A}

\define\fg{\frak g}

\define\fE{\frak E}

\define\tq{\ti q}

\define\tE{\ti E}
\define\tF{\ti F}

\define\sha{\sharp}

\define\CHE{Ch}
\define\KAS{Ka}
\define\QG{L1}
\define\ROOTS{L2}
\define\CAN{L3}
\define\QU{L4}

\head Introduction\endhead
\subhead 0.1\endsubhead
According to Drinfeld and Jimbo, the universal enveloping algebra of a
simple split Lie algebra $\fg$ over $\QQ$ admits a remarkable deformation 
$\UU$ (as a Hopf algebra over $\QQ(v)$, where $v$ is an indeterminate) called 
a quantized enveloping algebra. Moreover, the irreducible
finite dimensional $\fg$-modules admit quantum deformation to become simple
$\UU$-modules. In \cite{\CAN}, I found that these quantum deformations admit
canonical bases with very favourable properties (at least when $\fg$ is of type
$A,D$ or $E$) which give also rise by specialization to canonical bases of the
corresponding simple $\fg$-modules. (Later, Kashiwara \cite{\KAS} found 
another approach to the canonical bases.) In this paper we are interested in 
the canonical basis of the quantum deformation $\L$ of the adjoint 
representation of $\fg$. Before the introduction of the canonical bases, in 
\cite{\QG}, \cite{\ROOTS}, I found a basis of $\L$ in which the generators 
$E_i,F_i$ of $\UU$ act through matrices whose entries are polynomials in 
$\NN[v]$. By specialization, this gives rise to a basis of the adjoint 
representation of $\fg$ in which the Chevalley generators $e_i,f_i$ of $\fg$ 
act through matrices whose entries are natural numbers, in contrast with the 
more traditional treatments where a multitude of signs appear.

In this paper (Section 1) I will prove that the basis of $\L$ from \cite{\QG}, 
\cite{\ROOTS} coincides with the canonical basis of $\L$. I thank Meinolf Geck
for suggesting that I should write down this proof. As an application (Section
2), I will give a definition of the Chevalley group over a field $k$ 
associated to $\fg$ which seems to be simpler than Chevalley's original 
definition \cite{\CHE}.

\subhead 0.2\endsubhead
Let $I$ be a finite set with a given $\ZZ$-valued symmetric bilinear form 
$y,y'\m y\cdot y'$ on $Y=\ZZ[I]$ such that the symmetric matrix 
$(i\cdot j)_{i,j\in I}$ is positive definite and such that 
$i\cdot i/2\in\{1,2,3,\do\}$ for all $i\in I$, $i\cdot i/2=1$ for some 
$i\in I$ and $2\fra{i\cdot j}{i\cdot i}\in\{0,-1,-2,\do\}$ for all $i,j\in I$.
In the terminology of \cite{\QU, 1.1.1, 2.1.3}, this is a {\it Cartan datum} 
of finite type. We shall assume that our Cartan datum is irreducible (see 
\cite{\QU, 2.1.3}). Let $e$ be the maximum value of $i\cdot i/2$ for $i\in I$.
We have $e\in\{1,2,3\}$. Let $I^1=\{i\in I;i\cdot i/2=1\}$,
$I^e=\{i\in I;i\cdot i/2=e\}$. If $e=1$ we have clearly $I^1=I^e=I$; if $e>1$,
we have $I=I^1\sqc I^e$.

Let $X=\Hom(Y,\ZZ)$ and let $\la,\ra:Y\T X@>>>\ZZ$ be the obvious pairing. For
$j\in I$ we define $j'\in X$ by $\la i,j'\ra=2\fra{i\cdot j}{i\cdot i}$ for 
all $i\in I$. Let $v$ be an indeterminate. For $i\in I$ we set 
$v_i=v^{i\cdot i/2}$; for $n\in\ZZ$ we set 
$[n]_i=\fra{v_i^n-v_i^{-n}}{v_i-v_i\i}$; for $n\in\NN$ we set 
$[n]_i^!=\prod_{s=1}^n[s]_i$.

Note that when $i\in I^1$ we have $v_i=v$ and we write $[n]$ instead of 
$[n]_i$.

\subhead 0.3\endsubhead
Following Drinfeld and Jimbo we define $\UU$ to be the associative 
$\QQ(v)$-algebra with generators $E_i,F_i$ ($i\in I$), $K_y$ ($y\in Y$) and 
relations
$$K_yK_{y'}=K_{y+y'}\text{ for }y,y'\text{ in }Y,$$
$$K_iE_j=v^{\la i,j'\ra}E_jK_i\text{ for }i,j\text{ in }I,$$
$$K_iF_j=v^{-\la i,j'\ra}F_jK_i\text{ for }i,j\text{ in }I,$$
$$E_iF_j-F_jE_i=\d_{ij}\fra{K_i^{i\cdot i/2}-K_i^{-i\cdot i/2}}{v_i-v_i\i},$$
$$\sum_{p,p'\in\NN;p+p'=1-\la i,j'\ra}(-1)^{p'}\fra{[p+p']^!_i}
{[p]^!_i[p']^!_i}E_i^pE_jE_i^{p'}=0\text{ for }i\ne j\text{ in }I,$$
$$\sum_{p,p'\in\NN;p+p'=1-\la i,j'\ra}(-1)^{p'}\fra{[p+p']^!_i}
{[p]^!_i[p']^!_i}F_i^pF_jF_i^{p'}=0\text{ for }i\ne j\text{ in }I.$$
For $i\in I,s\in\NN$ we set $E_i^{(s)}=([s]^!_i)\i E_i^s$,
$F_i^{(s)}=([s]^!_i)\i F_i^s$.

By \cite{\QU, 3.1.12}, there is a unique $\QQ$-algebra isomorphism 
$\bar{}:\UU@>>>\UU$ such that $\bar E_i=E_i,\bar F_i=F_i$ for $i\in I$, 
$\bar K_y=K_{-y}$ for $y\in Y$ and $\ov{v^nu}=v^{-n}\bar u$ for all $u\in\UU$,
$n\in\ZZ$. 

\subhead 0.4\endsubhead
Let $W$ be the (finite) subgroup of $Aut(X)$ generated by the involutions
$s_i:\l\m\l-\la i,\l\ra i'$ of $X$ ($i\in I$). Let $R$ be the smallest 
$W$-stable subset of $X$ that contains $\{i';i\in I\}$. This is a finite set.
Let $R^+=\{\a\in R;\a\in\sum_i\NN i'\}$, $R^-=-R^+$. Let $R^1$ (resp. 
$R^e$ be the smallest $W$-stable subset of $X$ that contains $I^1$ (resp. 
$I^e$). Then $R^1,R^e$ are $W$-orbits. If $e=1$ we have $R=R^1=R^e$; if $e>1$ 
we have $R=R^1\sqc R^e$.

For $i\in I$ and $\a\in R$ let $p_{i,\a}$ be the largest integer $\ge0$ such 
that $\a,\a+i',\a+2i',\do,\a+p_{i,\a}i'$ belong to $R$ and let $q_{i,\a}$ be 
the largest integer $\ge0$ such that $\a,\a-i',\a-2i',\do,\a-q_{i,\a}i'$ 
belong to $R$. Then:

(a) $\la i,\a\ra=q_{i,\a}-p_{i,\a}$ and $p_{i,\a}+q_{i,\a}\le3$. 

(b) {\it If $p_{i,\a}+q_{i,\a}>1$, then we must have $p_{i,\a}+q_{i,\a}=e$, 
$i\in I^1$; moreover, $\a-q_{i,\a}i'\in R^e$, $\a+p_{i,\a}i'\in R^e$ and 
$\a+ki'\in R^1$ for $-q_{i,\a}<k<p_{i,\a}$.}

(c) {\it If $p_{i,\a}+q_{i,\a}=1$, then either both $\a-q_{i,\a}i'$, 
$\a+p_{i,\a}i'$ belong to $R^e$ or both belong to $R^1$.}
\nl
We define $h:R^+@>>>\NN$ by $h(\a)=\sum_{i\in I}n_i$ where
$\a=\sum_{i\in I}n_ii'$ with $n_i\in\NN$. There is a unique $\a_0\in R^+$ such
that $h(\a_0)$ is maximum. We then have $p_{i,\a_0}=0$ for all $i\in I$; it 
follows that $\la i,\a_0\ra\ge0$ for any $i\in I$. We have $\a_0\in R^e$.

\subhead 0.5\endsubhead
The $\UU$-module $\L:=\L_{\a_0}$ (see \cite{\QU, 3.5.6}) is well defined; it 
is simple, see \cite{\QU, 6.2.3}, and finite dimensional, see 
\cite{\QU, 6.3.4}. Let $\et=\et_{\a_0}\in\L$ be as in \cite{\QU, 3.5.7}. We 
have a direct sum decomposition (as a vector space) $\L=\op_{\l\in X}\L^\l$
where $\L^\l=\{x\in\L;K_yx=v^{\la y,\l\ra}x\qua\frl y\in Y\}$. Note that for 
$i\in I,\l\in X$ we have $E_iX^\l\sub X^{\l+i'}$, $F_iX^\l\sub X^{\l-i'}$. 
Moreover, we have $\dim\L^\a=1$ if $\a\in R$, $\dim\L^0=\sha(I)$ and 
$\L^\l=0$ if $\l\n R\cup\{0\}$.

Let $\BB$ be the canonical basis of $\L$ defined in \cite{\QU, 14.4.11}. We 
now state the following result in which $||$ denotes absolute value.

\proclaim{Theorem 0.6} (a) $\L$ has a unique $\QQ(v)$-basis
$\fE=\{X_\a;\a\in R\}\sqc\{t_i;i\in I\}$ such that (i)-(iii) below hold.

(i) $X_{\a_0}=\et$;

(ii) for $\a\in R$ we have $X_\a\in\L^\a$; for $i\in I$ we have $t_i\in\L^0$;

(iii) for any $i\in I$ the linear maps $E_i:\L@>>>\L$, $F_i:\L@>>>\L$, are 
given by
$$E_iX_\a=[q_{i,\a}+1]_iX_{\a+i'}\text{ if }\a\in R,p_{i,\a}>0,$$
$$E_iX_{-i'}=t_i,$$
$$E_iX_\a=0\text{ if }\a\in R,p_{i,\a}=0,\a\ne-i',$$
$$E_it_j=[|\la j,i'\ra|]_jX_{i'},\text{ if }j\in I,$$
$$F_iX_\a=[p_{i,\a}+1]_iX_{\a-i'}\text{ if }\a\in R,q_{i,\a}>0,$$
$$F_iX_{i'}=t_i,$$
$$F_iX_\a=0\text{ if }\a\in R,q_{i,\a}=0,\a\ne i',$$
$$F_it_j=[|\la j,i'\ra|]_jX_{-i'}\text{ if }j\in I.$$
(b) We have $\fE=\BB$.
\endproclaim
Note that the uniqueness of $\fE$ in (a) is straightforward. The existence of
$\fE$ is proved in \cite{\QG} under the assumption that $e=1$ and is stated in
\cite{\ROOTS} without assumption on $e$. We shall note use these results here.
Instead, in 1.15 we shall give a new proof (based on results in \cite{\QU}) of
the existence of $\fE$ at the same time as proving (b). 

\head 1. Proof of Theorem 0.6\endhead
\subhead 1.1\endsubhead
For any $\l\in X$, $\BB\cap\L^\l$ is a basis of $\L^\l$. In particular, for 
any $\a\in R$, $\BB\cap\L^\a$ is a single element; we denote it by $b^\a$. 

Let $\ca=\ZZ[v,v\i]$ and let $\L_\ca$ be the $\ca$-submodule of $\L$ generated
by $\BB$. It is known that $L_\ca$ is stable under $E_i^{(s)},F_i^{(s)}$ for
$i\in I,s\in N$.

By \cite{\QU, 19.3.4}, there is a unique $\QQ$-linear isomorphism 
$\bar{}:\L@>>>\L$ such that $\ov{u\et}=\bar u\et$ for all $u\in\UU$. By 
\cite{\QU, 19.1.2}, there is a unique bilinear form $(,):\L\T\L@>>>\QQ(v)$
such that $(\et,\et)=1$ and $(E_ix,x')=(x,v_iK_i^{i\cdot i/2}F_ix')$,
$(F_ix,x')=(x,v_iK_i^{-i\cdot i/2}E_ix')$, $(K_yx,x')=(x,K_yx')$ for all 
$i\in I,y\in Y$ and $x,x'$ in $\L$.

\subhead 1.2\endsubhead
By \cite{\QU, 19.3.5},

(a) {\it an element $b\in\L$ satisfies $\pm b\in\BB$ if and only if 
$b\in\L_\ca$, $\bar b=b$ and $(b,b)\in1+v\i\ZZ[v\i]$.}

\subhead 1.3\endsubhead
By \cite{\KAS} (see also \cite{\QU, 16.1.4}), 
for any $i\in I$ there is a unique $\QQ(v)$-linear map 
$\tF_i:\L@>>>\L$ such that the following holds: if $x\in\L^\l$, $E_ix=0$ and 
$s\in\NN$, then $\tF_i(F_i^{(s)}x)=F_i^{(s+1)}x$. Moreover, there is a unique 
$\QQ(v)$-linear map $\tE_i:\L@>>>\L$ such that the following holds: if 
$x\in\L^\l$, $F_ix=0$ and $s\in\NN$, then $\tE_i(E_i^{(s)}x)=E_i^{(s+1)}x$. 
Let $\AA=\QQ(v)\cap\QQ[[v\i]]$. Let $\L_\AA$ be the $\AA$-submodule of $\L$
generated by $\BB$. For any $x\in\L_\AA$ let $\un x$ be the image of $x$ in 
$\un\L:=\L_\AA/v\i\L_\AA$. Note that $\{\un b;b\in\BB\}$ is a $\QQ$-basis of 
$\un\L$. By \cite{\KAS} (see also \cite{\QU, 20.1.4}), 
for any $i\in I$, $\tF_i,\tE_i$ preserve 
$\L_\AA,v\i\L_\AA$ hence they induce $\QQ$-linear maps $\un\L@>>>\un\L$ 
(denoted again by $\tF_i,\tE_i$). From \cite{\KAS} (see also 
\cite{\QU, 20.1.4}) we see also that

(a) {\it $\tF_i:\un\L@>>>\un\L,\tE_i:\un\L@>>>\un\L$ act in the basis 
$\{\un b;b\in\BB\}$ by matrices with all entries in $\{0,1\}$.}
\nl
In the case where $e=1$, the results in this subsection are not needed; in
this case, instead of (a), we could use the positivity of the matrix entries
of $E_i:\L@>>>\L$, $F_i:\L@>>>\L$ proved in \cite{\QU, 22.1.7}.

\subhead 1.4\endsubhead
Let $\a\in R$, $i\in I$ be such that $q_{i,\a}=0,p=p_{i,\a}\ge1$. Then we have
$\la i,\a\ra=-p$. Let $Z^0=b^\a\in\L^\a$. We have $F_iZ^0\in\L^{\a-i'}$ hence 
$F_iZ^0=0$. We define $Z^k\in\L^{\a+ki'}$ for $k=1,\do,p$ by the inductive 
formula 

(a) $Z^k=[k]_i\i E_iZ^{k-1}=\tE_i^kZ^0$.
\nl
Using $F_iZ^0=0$ together with (a) and the commutation formula between 
$E_i,F_i$ we see by induction on $k$ that for $k=1,\do,p$ we have 

(b) $F_iZ^k=[p-k+1]_iZ^{k-1}$.

\subhead 1.5\endsubhead
We preserve the setup of 1.4. We show that for $k\in[0,p-1]$ we have
$$(Z^{k+1},Z^{k+1})=\fra{1-v_i^{-2p+2k}}{1-v_i^{-2k-2}}(Z^k,Z^k).\tag a$$
We have $E_iZ^k=[k+1]_iZ^{k+1}$ hence using 1.4(b):
$$\align&[k+1]_i^2(Z^{k+1},Z^{k+1})=(E_iZ^k,E_iZ^k)
=(Z^k,v_iK_i^{i\cdot i/2}F_iE_iZ^k)\\&=(Z^k,v_iK_i^{i\cdot i/2}E_iF_iZ^k)-
(Z^k,v_iK_i^{i\cdot i/2}\fra{K_i^{i\cdot i/2}-K_i^{-i\cdot i/2}}{v_i-v_i\i}Z^k)
\\&=(v_i^{\la i,\a+ki'\ra+1}[k]_i[p-k+1]_i
-\fra{v_i^{2\la i,\a+ki'\ra+1}-v_i}{v_i-v_i\i})(Z^k,Z^k)\\&
=(v_i^{-p+2k+1}[k]_i[p-k+1]_i-\fra{v_i^{-2p+4k+1}-v_i}{v_i-v_i\i})(Z^k,Z^k).
\endalign$$
We have
$$\align&(v_i-v_i\i)^2(v_i^{-p+2k+1}[k]_i[p-k+1]_i-\fra{v_i^{-2p+4k+1}-v_i}
{v_i-v_i\i})\\&=v_i^{-p+2k+1}(v_i^k-v_i^{-k})(v_i^{p-k+1}-v_i^{-p+k-1})
-(v_i^{-2p+4k+1}-v_i)(v_i-v_i\i)\\&
=v_i^{2k+2}-v_i^2-v_i^{-2p+4k}+v_i^{-2p+2k}-v_i^{-2p+4k+2}+v_i^2+v_i^{-2p+4k}-1
\\&=v_i^{2k+2}+v_i^{-2p+2k}-v_i^{-2p+4k+2}-1=(v_i^{-2p+2k}-1)(1-v_i^{2k+2}).
\endalign$$
Thus
$$(Z^{k+1},Z^{k+1})=
\fra{(v_i^{-2p+2k}-1)(1-v_i^{2k+2})}{(v_i^{k+1}-v_i^{-k-1})^2}(Z^k,Z^k)$$
and (a) follows.

\subhead 1.6\endsubhead
We preserve the setup of 1.4. We must have $p\in\{1,2,3\}$.

Assume first that $p=1$. From 1.5(a) we have $(Z^1,Z^1)=(Z^0,Z^0)$.

Assume now that $p=2$. Then from 0.4(b) we have $v_i=v$ and from 1.5(a) we 
have 

$(Z^1,Z^1)=\fra{1-v^{-4}}{1-v^{-2}}(Z^0,Z^0)$,

$(Z^2,Z^2)=\fra{1-v^{-2}}{1-v^{-4}}(Z^1,Z^1)=(Z^0,Z^0)$.
\nl
Assume next that $p=3$. Then from 0.4(b) we have $v_i=v$ and from 1.5(a) we 
have 

$(Z^1,Z^1)=\fra{1-v^{-6}}{1-v^{-2}}(Z^0,Z^0)$,

$(Z^2,Z^2)=(Z^1,Z^1)$.

$(Z^3,Z^3)=\fra{1-v^{-2}}{1-v^{-6}}(Z^2,Z^2)=(Z^0,Z^0)$.

\subhead 1.7\endsubhead
We preserve the setup of 1.6. We show:

(a) We have $Z^k=b^{\a+ki'}$ for $k=0,1,\do,p$.
\nl
Since $Z^0\in\BB$, we have
$Z^0\in\L_\ca$, $\bar Z^0=Z^0$, $(Z^0,Z^0)\in1+v\i\ZZ(v\i)$. From the 
formulas in 1.6 we see that $(Z^k,Z^k)\in1+v\i\ZZ(v\i)$ for $k=0,1,\do,p$. For
$k=1,\do,p$ we have $E_iZ^{k-1}=[k]_iZ^k$ hence for $k=0,1,\do,p$ we have 
$Z^k=E_i^{(k)}Z^0\in\L_\ca$. From $Z^k=E_i^{(k)}Z^0$ we see also that
$\bar Z^k=\ov{E_i^{(k)}}\ov{Z^0}=E_i^{(k)}Z^0=Z^k$. Using 1.2(a) we see that
$\e Z^k\in\BB$ for some $\e\in\{1,-1\}$. By 1.4(a), we have 
$\un{Z^k}=\tE_i^k\un{Z^0}$. Using this together with and 1.3(a), we see that 
$\e=1$ so that $Z^k\in\BB$. Since $Z^k\in\L^{\a+ki'}$, we see that 
$Z^k=b^{\a+ki'}$.

\subhead 1.8\endsubhead
Let $i\in I,\ti\a\in R$ be such that $p_{i,\ti\a}>0$ (or equivalently such that
$\ti\a+i'\in R$). We show:
$$E_ib^{\ti\a}=[q_{i,\ti\a}+1]_ib^{\ti\a+i'}\tag a$$
Let $\a=\ti\a-q_{i,\ti\a}i'\in R$. We have $q_{i,\a}=0$,
$p_{i,\a}=p_{i,\ti\a}+q_{i,\ti\a}>0$. We set $Z^0=b^\a$. We then define $Z^k$ 
with $k\in[1,p_{i,\a}]$ in terms of $\a,Z^0$ as in 1.4. Note that 
$E_iZ^{k-1}=[k]_iZ^k$ for any $k\in[1,p_{i,\a}]$. Taking $k=q_{i,\ti\a}+1$ (so
that $k\in[1,p_{i,\a}]$) we deduce
$$E_iZ^{q_{i,\ti\a}}=[q_{i,\ti\a}+1]_iZ^{q_{i,\ti\a}+1}.$$
By 1.7(a) we have $Z^{q_{i,\ti\a}}=b^{\ti\a}$, 
$Z^{q_{i,\ti\a}+1}=b^{\ti\a+i'}$. This proves (a).

\mpb

Here is a special case of (a); we assume that $i\ne j$ in $I$:

(b) If $\la j,i'\ra<0$ then $E_jb^{i'}=b^{i'+j'}$; if $\la j,i'\ra=0$ then 
$E_jb^{i'}=0$.
\nl
It is enough to use that $p_{j,i'}=-\la j,i'\ra$ (we have $q_{j,i'}=0$ since 
$i'-j'\n R$).

\subhead 1.9\endsubhead
Let $i\in I,\ti\a\in R$ be such that $q_{i,\ti\a}>0$ (or equivalently such that
$\ti\a-i'\in R$). We show:
$$F_ib^{\ti\a}=[p_{i,\ti\a}+1]_ib^{\ti\a-i'}.\tag a$$
Let $\a=\ti\a-q_{i,\ti\a}i'\in R$. We have $q_{i,\a}=0$,
$p_{i,\a}=p_{i,\ti\a}+q_{i,\ti\a}>0$. We set $Z^0=b^\a$. We then define $Z^k$ 
with $k\in[1,p_{i,\a}]$ in terms of $\a,Z^0$ as in 1.4.
Note that $F_iZ^k=[p_{i,\a}-k+1]_iZ^{k-1}$ for $k\in[1,p_{i,\a}]$. Taking
$k=q_{i,\ti\a}$ (so that $k\in[1,p_{i,\a}]$) we deduce
$$F_iZ^{q_{i,\ti\a}}=[p_{i,\ti\a}+1]_iZ^{q_{i,\ti\a}-1}.$$
By 1.7(a) we have $Z^{q_{i,\ti\a}}=b^{\ti\a}$, 
$Z^{q_{i,\ti\a}-1}=b^{\ti\a-i'}$. This proves (a).

\mpb

Here is a special case of (a); we assume that $i\ne j$ in $I$:

(b) If $\la j,i'\ra<0$ then $F_jb^{-i'}=b^{-i'-j'}$; if $\la j,i'\ra=0$, 
then $F_jb^{-i'}=0$.
\nl
It is enough to use that $q_{j,-i'}=\la j,-i'\ra$ (we have $p_{j,-i'}=0$ since
$-i'+j'\n R$).

\subhead 1.10\endsubhead
Let $i\in I$; we set $t_i=E_ib^{-i'}\in\L^0$. We show
$$F_it_i=(v_i+v_i\i)b^{-i'}.\tag a$$
Indeed,
$$\align&F_it_i=F_iE_ib^{-i'}=E_iF_ib^{-i'}-
\fra{K_i^{i\cdot i/2}-K_i^{-i\cdot i/2}}{v_i-v_i\i}b^{-i'}\\&=
\fra{v_i^2-v_i^{-2}}{v_i-v_i\i}b^{i'}=(v_i+v_i\i)b^{-i'}.\endalign$$
We show:
$$(t_i,t_i)=(1+v_i^{-2})(b^{-i'},b^{-i'}).\tag b$$
Indeed, using (a) we have
$$\align&(t_i,t_i)=(E_ib^{-i'},t_i)=(b^{-i'},v_iK_i^{i\cdot i/2}F_it_i)
=(b^{-i'},v_iK_i^{i\cdot i/2}(v_i+v_i\i)b^{-i'})=\\&
(v_i+v_i\i)v_i^{-\la i,i'\ra+1}(b^{-i'},b^{-i'})
=(1+v_i^{-2})(b^{-i'},b^{-i'}).\endalign$$
From (b) we see that $(t_i,t_i)\in 1+v\i\ZZ[v\i]$; from the definitions we 
have also $t_i\in\L_\ca$ and $\bar t_i=t_i$; it follows that $\e t_i\in\BB$
for some $\e\in\{1,-1\}$. Now from $t_i=E_ib^{-i'}$ and $F_ib^{-i'}=0$ we see 
that $t_i=\tE_ib^{-i'}$ hence $\un{t_i}=\tE_i\un{b^{-i'}}$. Using this 
together with 1.3(a) and we see that $\e=1$ hence
$$t_i\in\BB.\tag c$$
We show:
$$\text {If $i\ne j$, then }F_it_j=[-\la j,i'\ra]_jb^{-i'}.\tag d$$
We have $F_it_j=F_iE_jb^{-j'}=E_jF_ib^{-j'}$. This is $0$ if $\la i,j'\ra=0$ 
since by 1.9(b) we have $F_ib^{-j'}=0$ (so in this case (a) holds). Now assume 
that $\la i,j'\ra<0$. Then using 1.9(b) and 1.8(a) we have
$$E_jF_ib^{-j'}=E_jb^{-i'-j'}=[q_{j,-i'-j'}+1]_jb^{-i'}.$$
Note that $p_{j,-i'-j'}=1$ since $-i'-j'+j'\in R$, $-i'-j'+2j'\n R$.
Hence $q_{j,-i'-j'}-1=\la j,-i'-j'\ra=-2-\la j,i'\ra$ that is,
$q_{i,-i'-j'}+1=-\la j,i'\ra$. This completes the proof of (d).

We show:
$$(E_it_i,E_it_i)=[2]_i^2(b^{-i'},b^{-i'}).\tag e$$
Indeed, using (b) we have
$$\align&(E_it_i,E_it_i)=(t_i,v_iK_i^{i\cdot i/2}F_iE_it_i)
=(t_i,v_iK_i^{i\cdot i/2}E_iF_it_i)\\&-
(t_i,v_iK_i^{i\cdot i/2}\fra{K_i^{i\cdot i/2}-K_i^{-i\cdot i/2}}{v_i-v_i\i}t_i)
=[2]_i(t_i,v_iK_i^{i\cdot i/2}E_ib^{-i'})=[2]_i(t_i,v_iK_i^{i\cdot i/2}t_i)
\\&=[2]_i(t_i,v_it_i)=[2]_i^2(b^{-i'},b^{-i'}),\endalign$$
proving (e). 

From (e) we get $([2]_i\i E_it_i,[2]_i\i E_i t_i)\in 1+v\i\ZZ[v\i]$. We have 
$[2]_i\i E_it_i=E_i^{(2)}b^{-i'}\in\L_\ca$. Moreover, we have clearly 
$\ov{[2]_i\i E_it_i}=[2]_i\i E_it_i$. Using 1.2(a) we deduce that 
$\e[2]_i\i E_it_i\in\BB$ for some $\e\in\{1,-1\}$. Since 
$[2]_i\i E_it_i\in\L^{i'}$, we must have 
$\e[2]_i\i E_it_i=b^{i'}$. Thus we have
$\e E_i^{(2)}b^{-i'}=b^{i'}$. Since $F_ib^{-i'}=0$ it follows that
$\tE_i^2b^{-i'}=\e b^{i'}$ and $\tE_i^2\un{b^{-i'}}=\e\un{b^{i'}}$. Using 
1.3(a), we deduce that $\e=1$. Thus,
$$E_it_i=[2]_ib^{i'}.\tag f$$

\subhead 1.11\endsubhead
Let $i\in I$. We set $\ti t_i=F_ib^{i'}\in\L^0$. We show:
$$E_i\ti t_i=[2]_ib^{i'}.\tag a$$
Indeed,
$$E_i\ti t_i=E_iF_ib^{i'}=F_iE_ib^{i'}+
\fra{K_i^{i\cdot i/2}-K_i^{-i\cdot i/2}}{v_i-v_i\i}b^{i'}=
\fra{v_i^2-v_i^{-2}}{v_i-v_i\i}b^{i'}=[2]_ib^{i'}.$$
We show:
$$(\ti t_i,\ti t_i)=[2]_iv_i\i(b^{i'},b^{i'}).\tag b$$
Indeed, using (a) we have:
$$\align&(\ti t_i,\ti t_i)=(F_ib^{i'},\ti t_i)
=(b^{i'},v_iK_i^{-i\cdot i/2}E_i\ti t_i)
=(b^{i'},v_iK_i^{-i\cdot i/2}[2]_ib^{i'})\\&=
[2]_iv_i\i(b^{i'},b^{i'}).\endalign$$
From (b) we see that $(\ti t_i,\ti t_i)\in 1+v\i\ZZ[v\i]$; from the 
definitions we have also $\ti t_i\in\L_\ca$ and $\ov{\ti t}_i=\ti t_i$; using
1.2(a) we see that $\e\ti t_i\in\BB$ for some $\e\in\{1,-1\}$.
 From $\ti t_i=F_ib^{i'}$, $E_ib^{i'}=0$
we see that $\ti t_i=\tF_ib^{i'}$ hence $\un{\ti t_i}=\tF_i\un{b^{i'}}$. Using
this and 1.3(a) we deduce that $\e=1$ so that
$$\ti t_i\in\BB.\tag c$$
We show:
$$(\ti t_i,t_i)=\pm(1+v_i^{-2})(b^{i'},b^{i'}).\tag d$$
Indeed, using 1.10(f) we have
$$\align&(\ti t_i,t_i)=(F_ib^{i'},t_i)=(b^{i'},v_iK_i^{-i\cdot i/2}E_it_i)=
(b^{i'},v_iK_i^{-i\cdot i/2}[2]_ib^{i'})\\&=
v_i\i[2]_i(b^{i'},b^{i'})=(1+v_i^{-2})(b^{i'},b^{i'})\endalign$$
hence $(\ti t_i,t_i)\in 1+v\i\ZZ[v\i]$. 
If $\ti t_i\ne t_i$ then, since $\ti t_i\in\BB$ and $t_i\in\BB$, we would have
$(\ti t_i,t_i)\in v\i\ZZ[v\i]$ (see \cite{\QU, 19.3.3}), contradicting (d).
Thus we hve $\ti t_i=t_i$ and
$$F_ib^{i'}=t_i.\tag e$$ 
We show:
$$\text {If $i\ne j$, then }E_it_j=[-\la j,i'\ra]_jb^{i'}.\tag f$$
Using (e) we have $E_it_j=E_iF_jb^{j'}=F_jE_ib^{j'}$. This is $0$ if 
$\la i,j'\ra=0$ since by 1.8(b) we have $E_ib^{j'}=0$ (so in this case (f) 
holds). Now assume that $\la i,j'\ra<0$. Then using 1.8(b) and 1.9(a) we have
$$F_jE_ib^{j'}=F_jb^{i'+j'}=[p_{j,i'+j'}+1]_jb^{i'}.$$
Note that $q_{j,i'+j'}=1$ since $i'+j'-j'\in R$, $i'+j'-2j'\n R$. Hence 
$1-p_{j,i'+j'}=\la j,i'+j'\ra=2+\la j,i'\ra$ that is,
$p_{i,i'+j'}+1=-\la j,i'\ra$. This completes the proof of (f).

\subhead 1.12\endsubhead
We show:

(a) {\it If $\a\in R^1$, then 
$(b^\a,b^\a)=1+v^{-2}+\do+v^{-2(e-1)}=v^{-e+1}[e]$. If $\a\in R^e$, then 
$(b^\a,b^\a)=1$.}
\nl
Note that when $e=1$ we have $R^1=R^e$ and the two formulas in (a) are
compatible with each other.

We first prove (a) for $\a\in R^+$ by descending induction on $h(\a)$. If 
$h(\a)=h(\a_0)$ then $\a=\a_0$ and we have $b^\a=\et$ so that
$(b^\a,b^\a)=(\et,\et)=1$. Now assume that $\a\in R^+$, $h(a)<h(\a_0)$. We can
find $\a'\in R^+$, $i\in I$ such that $q_{i,\a'}=0$, $p=p_{i,\a'}\ge1$ and
$\a=\a'+ki'$ where $k\in\{0,1,\do,p-1\}$. Then $h(\a'+pi')>h(\a)$ hence 
$(\a'+pi',\a'+pi')$ is given by the formula in (a). Assume first that $p=1$. 
Then $\a=\a'$ and by 1.6 and 1.7(a) we have 
$(b^\a,b^\a)=(b^{\a'+i'},b^{\a'+i'})$. By 0.4(c), either both $\a,\a+i'$ 
belong to $R^e$ or both belong to $R^1$; (a) follows in this case. Next assume
that $p>1$. By 0.4(b) we have $p=e$ and $\a'+pi'\in R^e$.
Hence $(b^{\a'+pi'},b^{\a'+pi'})=1$.  
If $k=0$ then $\a\in R^e$ (see 0.4(b)) and by 1.6 and 1.7(a) we have
$(b^\a,b^\a)=(b^{\a'+pi'},b^{\a'+pi'})$; (a) follows in this case.
If $k>0$, $k<p$ then $\a\in R^1$ (see 0.4(b)) and by 1.6 and 1.7(a) we have
$(b^\a,b^\a)=(1+v^{-2}+\do+v^{-2(e-1)})(b^{\a'+pi'},b^{\a'+pi'})$;
(a) follows in this case.
This completes the proof of (a) assuming that $\a\in R^+$.

We now prove (a) for $\a\in R^-$ by induction on $h(-\a)\ge1$.
Let $i\in I$. Recall that $\ti t_i, t_i$ satisfy $\ti t_i=t_i$ (see 1.11),
$(t_i,t_i)=[2]_iv_i\i(b^{-i'},b^{-i'})$ (see 1.10(b)) and
$(\ti t_i,\ti t_i)=[2]_iv_i\i(b^{i'},b^{i'})$ (see 1.11(b)). It follows that
$$(b^{-i'},b^{-i'})=(b^{i'},b^{i'}).\tag b$$
In particular, (a) holds when $h(-\a)=1$.
We now assume that $\a\in R^-$ and $h(-\a)\ge2$. We can find $\a'\in R^-$, 
$i\in I$ such that $q_{i,\a'}=0$, $p=p_{i,\a'}\ge1$ and
$\a=\a'+ki'$ where $k\in\{0,1,\do,p-1\}$. Then 
$h(-(\a'+pi'))<h(-\a)$ hence 
$(\a'+pi',\a'+pi')$ is given by the formula in (a).
The rest of the proof is a repetition of the first part of the proof.
Assume first that $p=1$. Then $\a=\a'$ and by 1.6 and 1.7(a) we have 
$(b^\a,b^\a)=(b^{\a'+i'},b^{\a'+i'})$. By 0.4(c), either both $\a,\a+i'$ 
belong to $R^e$ or both belong to $R^1$; (a) follows in this case.
Next assume that $p>1$. By 0.4(b) we have $p=e$ and $\a'+pi'\in R^e$.
Hence $(b^{\a'+pi'},b^{\a'+pi'})=1$.  
If $k=0$ then $\a\in R^e$ (see 0.4(b)) and by 1.6 and 1.7(a) we have
$(b^\a,b^\a)=(b^{\a'+pi'},b^{\a'+pi'})$; (a) follows in this case.
If $k>0$, $k<p$ then $\a\in R^1$ (see 0.4(b)) and by 1.6 and 1.7(a) we have
$(b^\a,b^\a)=(1+v^{-2}+\do+v^{-2(e-1)})(b^{\a'+pi'},b^{\a'+pi'})$;
(a) follows in this case.
This completes the proof of (a) assuming that $\a\in R^-$; hence (a) is proved
in all cases.

\subhead 1.13\endsubhead
We show:

(a) {\it If $i\in I^1$ then $(t_i,t_i)=(1+v^{-2})(1+v^{-2}+\do+v^{-2(e-1)})$.
If $i\in I^e$ then $(t_i,t_i)=1+v_i^{-2}=1+v^{-2e}$.}
\nl
Note that when $e=1$ we have $I^1=I^e$ and the two formulas in (a) are
compatible with each other.

From 1.10(b) we have $(t_i,t_i)=[2]_iv_i\i(b^{-i'},b^{-i'})$. Using 1.12(a) we
see that (a) holds.

\mpb

In the remainder of this subsection we fix $i\ne j$ in $I$. We show:

(b) {\it If at least one of $i,j$ is in $I^1$ and $i\cdot j\ne0$ then
$(t_i,t_j)=v^{-e}[e]$. If both $i,j$ are in $I^e$ and $i\cdot j\ne0$ then
$(t_i,t_j)=v^{-e}$. If $i\cdot j=0$ then $(t_i,t_j)=0$.}
\nl
Using 1.10(d), we have
$$\align&(t_i,t_j)=(E_ib^{-i'},t_j)=(b^{-i'},v_iK_i^{i\cdot i/2}F_it_j)\\&
=[-\la j,i'\ra]_j(b^{-i'},v_iK_i^{i\cdot i/2}b^{-i'})
=v_i\i[-\la j,i'\ra]_j(b^{-i'},b^{-i'}).\endalign$$
We see that if $\la j,i'\ra=0$ then $(t_i,t_j)=0$.

Now assume that $\la j,i'\ra\ne0$. 

If $i\in I^e, j\in I^e$ then $\la j,i'\ra=-1$ and $(t_i,t_j)=v^{-e}$. 

If $i\in I^e, j\in I^1$ then $\la j,i'\ra=-e$ and $(t_i,t_j)=v^{-e}[e]$.

If $i\in I^1, j\in I^e$ then $(t_i,t_j)=(t_j,t_i)=v^{-e}[e]$.

If $i\in I^1, j\in I^1$ then $\la j,i'\ra=-1$ and 
$(t_i,t_j)=v\i(1+v^{-2}+\do+v^{-2(e-1)})=v^{-e}[e]$.
\nl
This completes the proof of (b).

\subhead 1.14\endsubhead
We show:

(a) {\it The elements $\{t_i;i\in I\}$ are distinct.}
\nl
Let $i\ne j$ in $I$. If we had $t_i=t_j$, then we would have 
$(t_i,t_j)\in1+v\i\ZZ[v\i]$, see 1.13(a). But 1.13(b) shows that 
$(t_i,t_j)\in v\i\ZZ[v\i]$. This completes the proof of (a).

Let $\fE=\{b^\a;\a\in R\}\sqc\{t_i,i\in I\}$. By (a), this is a subset of $\L$
rather than a multiset. We show:

(b) {\it We have $\BB=\fE$.}
\nl
Since $t_i\in\BB$ for any $i\in I$, we have $\fE\sub\BB$. Clearly we have
$\sha(\fE)=\sha(R)+\sha(I)$. Since we have also $\sha(\BB)=\sha(R)+\sha(I)$, 
it follows that $\fE=\BB$, proving (b).

\subhead 1.15\endsubhead
We prove the existence part of 0.6(a). It is enough to prove that the elements
$X_\a=b^\a$ and $t_i$ satisfy the requirements of 0.6(a). Now 0.6(a)(i) holds 
by definition; 0.6(a)(ii) is immediate; 0.6(a)(iii) has been verified earlier 
in this section. This proves the existence part of 0.6(a) and at the 
same time proves 0.6(b) (see 1.14(b)).

\head 2. Applications\endhead
\subhead 2.1\endsubhead
Let $i\in I$, $k\in\ZZ_{>0}$. From 0.6 we see that the action of $E_i^{(k)}$,
$F_i^{(k)}$ in the basis $\fE$ of $\L$ is given by the following formulas.
$$E_i^{(k)}X_\a=\fra{[q_{i,\a}+k]^!_i}{[q_{i,\a}]^!_i[k]^!_i}X_{\a+ki'}
\text{ if }\a\in R,\a\ne-i',k\le p_{i,\a},$$
$$E_i^{(k)}X_\a=0 \text{ if }\a\in R,\a\ne-i',k>p_{i,\a},$$
$$E_iX_{-i'}=t_i,E_i^{(2)}X_{-i'}=X_{i'}, E_i^{(k)}X_{-i'}=0\text{ if }k\ge3,$$
$$E_it_j=[|\la j,i'\ra|]_jX_{i'}, E_i^{(k)}t_j=0\text{ if }k\ge2,$$
$$F_i^{(k)}X_\a=\fra{[p_{i,\a}+k]^!_i}{[p_{i,\a}]^!_i[k]^!_i}X_{\a-ki'}
\text{ if }\a\in R,\a\ne i',k\le q_{i,\a},$$
$$F_i^{(k)}X_\a=0 \text{ if }\a\in R,\a\ne i',k>q_{i,\a},$$
$$F_iX_{i'}=t_i,F_i^{(2)}X_{i'}=X_{-i'}, F_i^{(k)}X_{i'}=0\text{ if }k\ge3,$$
$$F_it_j=[|\la j,i'\ra|]_jX_{-i'}, F_i^{(k)}t_j=0\text{ if }k\ge2.$$
In particular, we see that $E_i^{(k)}$, $F_i^{(k)}$ act through matrices with 
all entries in $\NN[v,v\i]$. (In the case where $e=1$ this is already known 
from \cite{\QU, 22.1.7}.)

\subhead 2.2\endsubhead
If $v$ is specialized to $1$, the $\UU$-module $\L$ becomes a simple
module over the universal enveloping algebra of a simple Lie algebra $\fg$
corresponding to the adjoint representation $\L|_{v=1}$ of $\fg$; this module
inherits a $\QQ$-basis $\{X_\a;\a\in R\}\sqc\{t_i;i\in I\}$ in which the 
elements $e_i,f_i$ of $\fg$ defined by $E_i,F_i$ act by matrices with entries 
in $\NN$. Let $z\in\QQ$. Then for $i\in I$, the exponentials 
$x_i(z)=\exp(ze_i),y_i(z)=\exp(zf_i)$ are well defined endomorphisms of 
$\L|_{v=1}$. Their action in the basis above can be described using the 
formulas in 2.1: 
$$x_i(z)X_\a=\sum_{0\le k\le p_{i,\a}}   
\fra{(q_{i,\a}+k)!}{q_{i,\a}!k!}z^kX_{\a+ki'}\text{ if }\a\in R,\a\ne-i',$$
$$x_i(z)X_{-i'}=X_{-i'}+zt_i+z^2X_{i'},$$
$$x_i(z)t_j=t_j+|\la j,i'\ra|zX_{i'}\text{ if }j\in I,$$ 
$$y_i(z)X_\a=\sum_{0\le k\le q_{i,\a}}
\fra{(p_{i,\a}+k)!}{p_{i,\a}!k!}z^kX_{\a-ki'}\text{ if }\a\in R,\a\ne i',$$
$$y_i(z)X_{i'}=X_{i'}+zt_i+z^2X_{-i'},$$
$$y_i(z)t_j=t_j+|\la j,i'\ra|zX_{-i'}\text{ if }j\in I.$$

\subhead 2.3\endsubhead
Now let $k$ be any field and let $V$ be the $k$-vector space with basis
$\{X_\a;\a\in R\}\sqc\{t_i;i\in I\}$. For any $i\in I$ and $z\in k$ we define
$x_i(z)\in GL(V)$, $y_i(z)\in GL(V)$ by the formulas in 2.2 (which involve 
only integer coefficients). The subgroup of $GL(V)$ generated by the elements 
$x_i(z),y_i(z)$ for various $i\in I,z\in k$ is the Chevalley group \cite{\CHE}
over $k$ associated to $\fg$.

\enddocument